\input amstex
\documentstyle{amsppt}
%----------------------------------------------------------------
% Title:     On operator fields in the bundle of Dirac spinors.
% Author:    Ruslan Sharipov
% Comments:  AmSTeX, 14 pages, amsppt style
% MSC-class: 53B30, 81T20, 81R25
%----------------------------------------------------------------
%           Replacement for output macro definition
%
\catcode`@=11
\redefine\output@{%
  \def\break{\penalty-\@M}\let\par\endgraf
  \ifodd\pageno\global\hoffset=105pt\else\global\hoffset=8pt\fi  
  \shipout\vbox{%
    \ifplain@
      \let\makeheadline\relax \let\makefootline\relax
    \else
      \iffirstpage@ \global\firstpage@false
        \let\rightheadline\frheadline
        \let\leftheadline\flheadline
      \else
        \ifrunheads@ %\let\makefootline\relax
        \else \let\makeheadline\relax
        \fi
      \fi
    \fi
    \makeheadline \pagebody \makefootline}%
  \advancepageno \ifnum\outputpenalty>-\@MM\else\dosupereject\fi
}
\catcode`\@=\active
%----------------------------------------------------------------
\nopagenumbers
\def\negskp{\hskip -2pt}
\def\tr{\operatorname{tr}}
\accentedsymbol\bd{\kern 2pt\bar{\kern -2pt d}}
\accentedsymbol\bbd{\kern 2pt\bar{\kern -2pt\bold d}}
\def\vtrule{\vrule height 12pt depth 6pt}
\def\vtttrule{\vrule height 12pt depth 19pt}
\def\boxit#1#2{\vcenter{\hsize=122pt\offinterlineskip\hrule
  \line{\vtttrule\hss\vtop{\hsize=120pt\centerline{#1}\vskip 5pt
  \centerline{#2}}\hss\vtttrule}\hrule}}
\def\blue#1{#1}

\catcode`#=11\def\diez{#}\catcode`#=6
\catcode`_=11\def\podcherkivanie{_}\catcode`_=8
%\catcode`~=11\def\volna{~}\catcode`~=\active
\def\mycite#1{\cite{\blue{#1}}\immediate\special{ps:
     ShrHPSdict begin /ShrBORDERthickness 0 def}}

\def\mytag#1{%
    \tag#1}
\def\mythetag#1{\thetag{\blue{#1}}\immediate\special{ps:
     ShrHPSdict begin /ShrBORDERthickness 0 def}}
\def\myrefno#1{\no#1}
\def\myhref#1#2{\blue{#2}\immediate\special{ps:
     ShrHPSdict begin /ShrBORDERthickness 0 def}}
\def\myEarXivlink{\myhref{http://arXiv.org}{http:/\negskp/arXiv.org}}

\def\mytheorem#1{\csname proclaim\endcsname{Theorem #1}}
\def\mythetheorem#1{\blue{#1}\immediate\special{ps:
     ShrHPSdict begin /ShrBORDERthickness 0 def}}
\def\mylemma#1{\csname proclaim\endcsname{Lemma #1}}

\def\mycorollary#1{\csname proclaim\endcsname{Corollary #1}}

\def\mydefinition#1{\definition{Definition #1}}

%----------------------------------------------------------------
% Cyrillic fonts definition
%\font\eightcyr=wncyr8
%----------------------------------------------------------------
\pagewidth{360pt}
\pageheight{606pt}
\topmatter
\title
On operator fields\\in the bundle of Dirac spinors.
\endtitle
\author
R.~A.~Sharipov
\endauthor
\address 5 Rabochaya street, 450003 Ufa, Russia\newline
\vphantom{a}\kern 12pt Cell Phone: +7(917)476 93 48
\endaddress
\email \vtop to 30pt{\hsize=280pt\noindent
\myhref{mailto:r-sharipov\@mail.ru}
{r-sharipov\@mail.ru}\newline
\myhref{mailto:R\podcherkivanie Sharipov\@ic.bashedu.ru}
{R\_\hskip 1pt Sharipov\@ic.bashedu.ru}\vss}
\endemail
\urladdr
\vtop to 20pt{\hsize=280pt\noindent
\myhref{http://www.geocities.com/r-sharipov}
{http:/\negskp/www.geocities.com/r-sharipov}\newline
\myhref{http://www.freetextbooks.boom.ru/index.html}
{http:/\negskp/www.freetextbooks.boom.ru/index.html}\newline
\myhref{http://sovlit2.narod.ru/}
{http:/\negskp/sovlit2.narod.ru}\vss}
\endurladdr
\abstract
    Operator fields in the bundle of Dirac spinors and their 
conversion to spatial fields are considered. Some commutator
equations are studied with the use of the conversion technique.
\endabstract
\endtopmatter
\loadbold
%\loadeufb
\TagsOnRight
\document

\head
1. Introduction.
\endhead
    The bundle of Dirac spinors is used for describing particles with 
half-integer spin in general relativity and in quantum field theory.
It is a special four-dimensional complex vector-bundle over the 
space-time manifold $M$. Let's remind that the space-time manifold $M$ 
itself is a four-dimensional real manifold equipped with a Minkowski 
type metric $\bold g$ of the signature $(+,-,-,-)$. Apart from $\bold g$, 
the space-time manifold $M$ is equipped with two other geometric structures
--- the orientation and the polarization. The orientation distinguishes
right quadruples of tangent vectors from left ones, while the polarization
distinguishes future and past half light cones in tangent spaces at each 
point of $M$.\par 
     The bundle of Dirac spinors is denoted $DM$. It is equipped with four
basic spin-tensorial fields in addition to $\bold g$. They are presented
in the following table.
$$
\vcenter{\hsize 10cm
\offinterlineskip\settabs\+\indent
\vtrule
\hskip 1.2cm &\vtrule % Quantity
\hskip 5.2cm &\vtrule % Unit
\hskip 2.8cm &\vtrule % Relation
\cr\hrule 
\+\vtrule
\hfill\,Symbol\hfill&\vtrule
\hfill Name\hfill &\vtrule
\hfill Spin-tensorial\hfill &\vtrule\cr
\vskip -0.2cm
\+\vtrule
\hfill &\vtrule
\hfill \hfill&\vtrule
\hfill type\hfill&\vtrule\cr\hrule
\+\vtrule
\hfill $\bold d$\hfill&\vtrule
\hfill Skew-symmetric metric tensor\hfill&\vtrule
\hfill $(0,2|0,0|0,0)$\hfill&\vtrule\cr\hrule
\+\vtrule
\hfill$\bold H$\hfill&\vtrule
\hfill Chirality operator\hfill&\vtrule
\hfill $(1,1|0,0|0,0)$\hfill&\vtrule\cr\hrule
\+\vtrule
\hfill$\bold D$\hfill&\vtrule
\hfill Dirac form\hfill&\vtrule
\hfill $(0,1|0,1|0,0)$\hfill&\vtrule\cr\hrule
\+\vtrule
\hfill$\boldsymbol\gamma$\hfill&\vtrule
\hfill Dirac $\gamma$-field\hfill&\vtrule
\hfill $(1,1|0,0|1,0)$\hfill&\vtrule\cr\hrule
}
$$
The metric tensor $\bold g$ itself is interpreted as a spin-tensorial 
field of the spin-tensorial type $(0,0|0,0|0,2)$.\par
    In this paper, saying an operator field, we assume a spin-tensorial
field $\bold F$ of the spin-tensorial type $(1,1|0,0|0,0)$. In the 
coordinate form it is presented by a matrix $F^a_b$, where $a$ and $b$
are two spinor indices. Each operator field $\bold F$ in the bundle of 
Dirac spinors has a unique presentation of the following form:
$$
\hskip -2em
\gathered
F^a_b=u\,\delta^a_b+v\,H^a_b+\sum^3_{k=0}\gamma^{ak}_b\,u_k\,+\\
+\sum^3_{k=0}\sum^4_{c=1}H^a_c\,\gamma^{ck}_b\,v_k+\sum^3_{p=0}
\sum^3_{q=0}\sum^4_{c=1}\gamma^{ap}_c\,\gamma^{cq}_b
\,w_{p\kern 0.4pt q}.
\endgathered
\mytag{1.1}
$$
Here $u$ and $v$ are two scalar fields, $u_k$ and $v_k$ are the 
components of two covectorial fields $\bold u$ and $\bold v$, and
$w_{p\kern 0.4pt q}$ are the components of a skew-symmetric 
tensorial field $\bold w$. Through $H^a_b$ in \mythetag{1.1} we
denote the components of the chirality operator $\bold H$, while 
$\gamma^{ak}_b$ are the components of the Dirac $\gamma$-filed
$\boldsymbol\gamma$.\par
    The formula \mythetag{1.1} is a conversion formula associating
the spin-tensorial operator field $\bold F$ with the purely tensorial
fields $u$, $v$, $\bold u$, $\bold v$, $\bold w$. The presentation
\mythetag{1.1} is well-known (see \S\,28 in \mycite{1}). The main
goal of this paper is to study some special commutator equations
for operator fields in terms of their associated tensorial fields.
\head
2. Frames and coordinate presentations\\
of the basic fields.
\endhead
\mydefinition{2.1} A {\it spatial frame\/} is a quadruple of vector fields
$\boldsymbol\Upsilon_0,\,\boldsymbol\Upsilon_1,\,\boldsymbol\Upsilon_2,
\,\boldsymbol\Upsilon_3$ defined in some open domain of the space-time
manifold $M$ and linearly independent at each point of their domain.
\enddefinition
\mydefinition{2.2} A spatial frame $\boldsymbol\Upsilon_0,\,\boldsymbol
\Upsilon_1,\,\boldsymbol\Upsilon_2,\,\boldsymbol\Upsilon_3$ is called
a {\it right frame\/} if at each point of its domain its vectors 
$\boldsymbol\Upsilon_0,\,\boldsymbol\Upsilon_1,\,\boldsymbol\Upsilon_2,
\,\boldsymbol\Upsilon_3$ form a right quadruple in the sense of the 
orientation in $M$.
\enddefinition
\mydefinition{2.3} A spatial frame $\boldsymbol\Upsilon_0,\,\boldsymbol
\Upsilon_1,\,\boldsymbol\Upsilon_2,\,\boldsymbol\Upsilon_3$ is called
an {\it orthonormal frame\/} if the metric tensor $\bold g$ is presented 
by the standard Minkowski matrix in this frame:
$$
\hskip -2em
g_{ij}=g^{ij}=\Vmatrix 1 & 0 & 0 & 0\\0 & -1 & 0 & 0\\
0 & 0 & -1 & 0\\0 & 0 & 0 & -1\endVmatrix.
\mytag{2.1}
$$
\enddefinition
\noindent
In physical literature the matrix \mythetag{2.1} is often denoted by
$\eta_{ij}$. However, this is not a good tradition. I prefer to use
the symbol $g$ for the components of the metric tensor irrespective
to the choice of an orthonormal or a non-orthonormal frame.
\mydefinition{2.4} An orthonormal spatial frame $\boldsymbol\Upsilon_0,
\,\boldsymbol\Upsilon_1,\,\boldsymbol\Upsilon_2,\,\boldsymbol\Upsilon_3$
is called {\it positively polarized} if its first vector $\boldsymbol
\Upsilon_0$ is a time-like vector directed to the future in the sense
of the polarization in $M$.
\enddefinition
     A positively polarized right orthonormal frame $\boldsymbol
\Upsilon_0,\,\boldsymbol\Upsilon_1,\,\boldsymbol\Upsilon_2,\,\boldsymbol
\Upsilon_3$ is a typical choice when dealing with spinors. Note, however,
that in general case of a non-flat space-time $M$ such a frame is not
holonomic, i\.\,e\. its vector fields do not commute:
$$
\hskip -2em
[\boldsymbol\Upsilon_{\!i},\boldsymbol\Upsilon_{\!j}]=\sum^3_{k=0}
c^{\,k}_{ij}\,\boldsymbol\Upsilon_{\!k}.
\mytag{2.2}
$$
The coefficients $c^{\,k}_{ij}$ in \mythetag{2.2} are called the
{\it commutation coefficients\/} of the frame $\boldsymbol\Upsilon_0,
\,\boldsymbol\Upsilon_1,\,\boldsymbol\Upsilon_2,\,\boldsymbol\Upsilon_3$.
This frame is called {\it holonomic} if all of its commutation coefficients
are identically zero. Otherwise, it is called a {\it non-holonomic} frame.
\par
    It is known that the metric $\bold g$ induces the $4$-form $\boldsymbol
\omega$ in $M$ which is called the {\it volume form} or the {\it volume
tensor}. This differential form is used for integration over $M$. In the 
coordinate form the volume tensor $\boldsymbol\omega$ is given by the formula
$$
\hskip -2em
\omega_{ijkm}=\pm\sqrt{-\det(\smash{g_{ij}})}\,\varepsilon_{ijkm},
\mytag{2.3}
$$
where $\varepsilon$ is the Levi-Civita symbol:
$$
\hskip -2em
\varepsilon_{ijkm}=\varepsilon^{ijkm}
=\cases\ \ 1 &\vtop{\baselineskip=2pt\lineskip=2pt\hsize=140pt\noindent if 
$(ijkm)$ is an even permutation of the numbers $0,\,1,\,2,\,3$;}\\
-1 &\vtop{\baselineskip=2pt\lineskip=2pt\hsize=140pt\noindent if $(ijkm)$
is an odd permutation of the numbers $0,\,1,\,2,\,3$;}\\
\ \ 0 &\text{in all other cases.}
\endcases
\mytag{2.4}
$$
Typically $\boldsymbol\omega$ is treated as a pseudotensor. However, we
assume $M$ to be an orientable manifold with a fixed orientation. In this
case we can fix the choice of sign in \mythetag{2.3} by setting plus for 
right frames and setting minus for left frames. Therefore, we treat 
$\boldsymbol\omega$ as a tensor.\par
     The dual volume tensor is denoted by the same symbol $\boldsymbol
\omega$. Its components are produced from $\omega_{ijkm}$ by means of
the standard index raising procedure:
$$
\hskip -2em
\omega^{ijkm}=\sum^3_{p=0}\sum^3_{q=0}\sum^3_{r=0}\sum^3_{s=0}
\omega_{pqrs}\,g^{pi}\,g^{qj}\,g^{rk}\,g^{sm}.
\mytag{2.5}
$$
Applying the formula \mythetag{2.3} to \mythetag{2.5}, we derive the 
formula
$$
\hskip -2em
\omega^{ijkm}=\mp\sqrt{-\det(g^{ij})}\,\varepsilon^{ijkm},
\mytag{2.6}
$$
where $\varepsilon$ again is the Levi-Civita symbol \mythetag{2.4}. In 
the case of a right orthonormal frame the formulas \mythetag{2.3} and 
\mythetag{2.6} are reduced to 
$$
\xalignat 2
&\hskip -2em
\omega_{ijkm}=\varepsilon_{ijkm},
&&\omega^{ijkm}=-\varepsilon^{ijkm}.
\mytag{2.7}
\endxalignat
$$\par
\mydefinition{2.5} A spinor frame in the bundle of Dirac spinors $DM$ is 
a quadruple $\boldsymbol\Psi_1,\,\boldsymbol\Psi_2,\,\boldsymbol\Psi_3,
\,\boldsymbol\Psi_4$ of smooth sections of this bundle over some open 
domain of $M$ linearly independent at each point of this domain.
\enddefinition
\mydefinition{2.6} A spinor frame $\boldsymbol\Psi_1,\,\boldsymbol\Psi_2,
\,\boldsymbol\Psi_3,\,\boldsymbol\Psi_4$ is called an {\it orthonormal 
frame\/} the spinor metric $\bold d$ is presented by the following matrices 
in this frame:
$$
\xalignat 2
&\hskip -2em
d_{ij}=\Vmatrix 0 & 1 & 0 & 0\\-1 & 0 & 0 & 0\\
0 & 0 & 0 & -1\\0 & 0 & 1 & 0\endVmatrix,
&&d^{ij}=\Vmatrix 0 & -1 & 0 & 0\\1 & 0 & 0 & 0\\
0 & 0 & 0 & 1\\0 & 0 & -1 & 0\endVmatrix.
\mytag{2.8}
\endxalignat
$$
\enddefinition
\noindent
The matrices \mythetag{2.8} are inverse to each other. They present the
spinor metric $\bold d$ and its dual metric in an orthonormal spinor frame.
\pagebreak Irrespective to the choice of a spinor frame (orthonormal or 
non-orthonormal) the components of the spinor metric $\bold d$ are used 
for lowering spinor indices. The components of the dual spinor metric are 
used for raising spinor indices.\par
\mydefinition{2.7} A spinor frame $\boldsymbol\Psi_1,\,\boldsymbol\Psi_2,
\,\boldsymbol\Psi_3,\,\boldsymbol\Psi_4$ of the bundle $DM$ is called 
a {\it chiral frame\/} if the chirality operator $\bold H$ given by the 
following matrix in this frame:
$$
\hskip -2em
H^{\kern 0.5pti}_{\kern -0.5pt j}=\Vmatrix 1 & 0 & 0 & 0\\0 & 1 & 0 & 0\\
0 & 0 & -1 & 0\\0 & 0 & 0 & -1\endVmatrix.
\mytag{2.9}
$$
\enddefinition
\mydefinition{2.8} A spinor frame $\boldsymbol\Psi_1,\,\boldsymbol\Psi_2,
\,\boldsymbol\Psi_3,\,\boldsymbol\Psi_4$ of the Dirac bundle $DM$ is
called a {\it self-adjoint frame\/} if the Hermitian metric tensor 
$\bold D$ (the Dirac form) is represented by the following matrix in this
frame:
$$
\hskip -2em
D_{i\bar j}=\Vmatrix 0 & 0 & 1 & 0\\0 & 0 & 0 & 1\\
1 & 0 & 0 & 0\\0 & 1 & 0 & 0\endVmatrix.
\mytag{2.10}
$$
\enddefinition
\mydefinition{2.9}{\it Canonically orthonormal chiral frames} in $DM$ are 
those which are orthonormal, chiral, and self-adjoint simultaneously .
\enddefinition
     Canonically orthonormal chiral frames in $DM$ do exist. Moreover, 
each such frame is canonically associated with some positively polarized
right orthonormal frame in $TM$. Apart from canonically orthonormal chiral 
frames, there are three other special types of frames in $DM$. All of these 
frame types and their associated frame types in $TM$ are listed in the
following diagram.
$$
\hskip -2em
\aligned
&\boxit{Canonically orthonormal}{chiral frames}\to
\boxit{Positively polarized}{right orthonormal frames}\\
&\boxit{$P$-reverse}{anti-chiral frames}\to
\boxit{Positively polarized}{left orthonormal frames}\\
&\boxit{$T$-reverse}{anti-chiral frames}\to
\boxit{Negatively polarized}{left orthonormal frames}\\
&\boxit{$PT$-reverse}{chiral frames}\to
\boxit{Negatively polarized}{right orthonormal frames}
\endaligned
\mytag{2.11}
$$
More details concerning the diagram \mythetag{2.11} can be found in 
\mycite{2}. In this paper we shall use canonically orthonormal chiral 
frames in $DM$ and their associated positively polarized right 
orthonormal frames in $TM$ only. They are sufficient for our 
purposes.\par
     The bundle of Dirac spinors $DM$ is a complex vector bundle. 
Therefore it is equipped with the involution of complex conjugation
$\tau$ that acts upon spin-tensorial fields and changes their
spin tensorial type as follows:
$$
\hskip -2em
\CD
@>\tau>>\\
\vspace{-4ex}
(\alpha,\beta|\nu,\gamma|r,s)@.(\nu,\gamma|\alpha,\beta|r,s).\\
\vspace{-4.2ex}
@<<\tau< 
\endCD
\mytag{2.12}
$$
As we see in \mythetag{2.12}, the involution $\tau$ exchanges spinor 
and conjugate spinor indices. In the coordinate form it acts through
complex conjugation upon the components of spin-tensors.\par
     Applying $\tau$ to $\bold H$ and $\bold d$ we get two other
basic fields $\bar\bold H=\tau(\bold H)$ and $\bbd=\tau(\bold d)$.
They are called the {\it conjugate chirality operator} and the
{\it conjugate spinor metric} respectively. In a canonically orthonormal 
chiral frame $\boldsymbol\Psi_1,\,\boldsymbol\Psi_2,\,\boldsymbol\Psi_3,
\,\boldsymbol\Psi_4$ the conjugate chirality operator $\bar\bold H$ is
given by the matrix 
$$
\hskip -2em
\bar H^{\kern 0.5pt\bar i}_{\kern -0.5pt\bar j}=\Vmatrix 1 & 0 & 0 & 0
\\0 & 1 & 0 & 0\\0 & 0 & -1 & 0\\0 & 0 & 0 & -1\endVmatrix.
\mytag{2.13}
$$
The conjugate spinor metric in such a frame is given by the matrices
$$
\xalignat 2
&\hskip -2em
\bd_{\bar i\bar j}=\Vmatrix 0 & 1 & 0 & 0\\-1 & 0 & 0 & 0\\
0 & 0 & 0 & -1\\0 & 0 & 1 & 0\endVmatrix,
&&\bd^{\bar i\bar j}=\Vmatrix 0 & -1 & 0 & 0\\1 & 0 & 0 & 0\\
0 & 0 & 0 & 1\\0 & 0 & -1 & 0\endVmatrix.
\mytag{2.14}
\endxalignat
$$
Though the matrix \mythetag{2.13} coincides with the matrix
\mythetag{2.9} and the matrices \mythetag{2.14} coincide with 
the matrices \mythetag{2.8}, $\bar\bold H\neq\bold H$ and $\bbd\neq
\bold d$ because the spin tensorial types of these fields are different.
The coincidence of their matrices occurring in a special frame is 
destroyed in an arbitrary non-special frame.\par
     The Dirac form $\bold D$ is invariant with respect to the involution
$\tau$, i\.\,e\. we have $\tau(\bold D)=\bold D$. In the coordinate form
this equality is written as follows:
$$
\hskip -2em
D_{i\bar j}=\overline{D_{\bar j\kern 0.5pt i}}
\mytag{2.15}
$$
The equality \mythetag{2.15} is easily derived from \mythetag{2.10}. Being
derived in a special frame, it remains valid in an arbitrary frame too.\par
     In order to present the Dirac $\gamma$-field in the coordinate form
we need to fix two frames --- some spinor frame $\boldsymbol\Psi_1,\,
\boldsymbol\Psi_2,\,\boldsymbol\Psi_3,\,\boldsymbol\Psi_4$ and some spatial 
frame $\boldsymbol\Upsilon_0,\,\boldsymbol\Upsilon_1,\,\boldsymbol\Upsilon_2$,
$\boldsymbol\Upsilon_3$. We do it according to the first line in the diagram 
\mythetag{2.11}. In other words, we choose some canonically orthonormal chiral 
frame $\boldsymbol\Psi_1,\,\boldsymbol\Psi_2,\,\boldsymbol\Psi_3,\,\boldsymbol
\Psi_4$ in $DM$ and take its associated positively polarized right orthonormal
frame in $TM$. Then the Dirac $\gamma$-field $\boldsymbol\gamma$ is given by
the following four matrices:
$$
\xalignat 2
&\hskip -2em
\gamma^0=\Vmatrix 0&0&1&0\\ 0&0&0&1\\ 1&0&0&0\\ 0&1&0&0\endVmatrix,
&&\gamma^1=\Vmatrix 0&0&0&-1\\ 0&0&-1&0\\ 0&1&0&0\\ 1&0&0&0\endVmatrix,\\
\vspace{-1.5ex}
&&&\mytag{2.16}\\
\vspace{-1.5ex}
&\hskip -2em
\gamma^2=\Vmatrix 0&0&0&i\\ 0&0&-i&0\\ 0&-i&0&0\\ i&0&0&0\endVmatrix,
&&\gamma^3=\Vmatrix 0&0&-1&0\\ 0&0&0&1\\ 1&0&0&0\\ 0&-1&0&0\endVmatrix.
\endxalignat
$$
The number of a matrix in \mythetag{2.16} is determined by the spatial
index $k$ of the component $\gamma^{ak}_b$. Two spinor indices $a$ and 
$b$ specify the position of the component $\gamma^{ak}_b$ within the
matrix $\gamma^k$.\par
\head
3. Some algebraic relationships for the basic fields.
\endhead
    First of all let's note that the square of the chirality operator 
$\bold H$ is equal to the unit operator $\bold 1$. The same is true for 
the conjugate chirality operator $\bar\bold H$:
$$
\xalignat 2
&\hskip -2em
\bold H^2=\bold 1,
&&\bar\bold H^2=\bold 1.
\mytag{3.1}
\endxalignat
$$
In the coordinate form the identities \mythetag{3.1} are written as 
$$
\xalignat 2
&\hskip -2em
\sum^4_{c=1}H^a_c\,H^c_b=\delta^a_b,
&&\sum^4_{\bar c=1}\bar H^{\bar a}_{\bar c}\,\bar H^{\bar c}_{\bar b}
=\delta^{\bar a}_{\bar b}.
\mytag{3.2}
\endxalignat
$$
The formulas \mythetag{3.2} immediately follow from \mythetag{2.9} and
\mythetag{2.13}.\par
     When some spatial basis $\boldsymbol\Upsilon_0,\,\boldsymbol
\Upsilon_1,\,\boldsymbol\Upsilon_2,\,\boldsymbol\Upsilon_3$ is fixed,
the matrices \mythetag{2.16} can be treated as the components of four 
operators $\boldsymbol\gamma^0,\,\boldsymbol\gamma^1,\,\boldsymbol
\gamma^2,\,\boldsymbol\gamma^3$ acting in fibers of the bundle $DM$. 
These operators satisfy the following well-known relationships:
$$
\hskip -2em
\{\boldsymbol\gamma^p,\boldsymbol\gamma^q\}=2\,g^{p\kern 0.5pt q}\,\bold 1.
\mytag{3.3}
$$
The curly brackets in \mythetag{3.3} denotes the anticommutator of operators.
In the coordinate form the relationships \mythetag{3.3} are written as
$$
\sum^4_{c=1}\gamma^{ap}_c\,\gamma^{cq}_b+\sum^4_{c=1}\gamma^{aq}_c
\,\gamma^{cp}_b=2\,g^{p\kern 0.5pt q}\,\delta^a_b.
$$
These relationships are proved by direct calculations with the
use of \mythetag{2.16}. Lowering the upper spatial index of $\boldsymbol
\gamma^0,\,\boldsymbol\gamma^1,\,\boldsymbol\gamma^2,\,\boldsymbol\gamma^3$ 
by means of the metric $\bold g$, we get other four operators $\boldsymbol
\gamma_0,\,\boldsymbol\gamma_1,\,\boldsymbol\gamma_2,\,\boldsymbol\gamma_3$:
$$
\hskip -2em
\boldsymbol\gamma_k=\sum^3_{q=0}g_{kq}\,\boldsymbol\gamma^q.
\mytag{3.4}
$$
The operators $\boldsymbol\gamma_0,\,\boldsymbol\gamma_1,\,\boldsymbol
\gamma_2,\,\boldsymbol\gamma_3$ satisfy the following relationships:
$$
\xalignat 2
&\hskip -2em
\{\boldsymbol\gamma^p,\boldsymbol\gamma_q\}=2\,\delta^p_q\,\bold 1,
&&\{\boldsymbol\gamma_p,\boldsymbol\gamma_q\}=2\,g_{p\kern 0.5pt q}
\,\bold 1.
\mytag{3.5}
\endxalignat 
$$
The relationships \mythetag{3.5} are easily derived from \mythetag{3.3}
with the use of \mythetag{3.4}. The explicit matrix presentation for the
operators \mythetag{3.4} is derived from \mythetag{2.16}:
$$
\xalignat 2
&\hskip -2em
\gamma_0=\Vmatrix 0&0&1&0\\ 0&0&0&1\\ 1&0&0&0\\ 0&1&0&0\endVmatrix,
&&\gamma_1=\Vmatrix 0&0&0&1\\ 0&0&1&0\\ 0&-1&0&0\\ -1&0&0&0\endVmatrix,\\
\vspace{-1.5ex}
&&&\mytag{3.6}\\
\vspace{-1.5ex}
&\hskip -2em
\gamma_2=\Vmatrix 0&0&0&-i\\ 0&0&i&0\\ 0&i&0&0\\ -i&0&0&0\endVmatrix,
&&\gamma_3=\Vmatrix 0&0&1&0\\ 0&0&0&-1\\ -1&0&0&0\\ 0&1&0&0\endVmatrix.
\endxalignat
$$
Like the formulas \mythetag{2.16}, thee formulas \mythetag{3.6} are valid
if we choose some canonically orthonormal chiral frame in $DM$ and its 
associated positively polarized right orthonormal frame in $TM$.
\par
      Note that the chirality operator $\bold H$ can be expressed through
the operators $\boldsymbol\gamma^0,\,\boldsymbol\gamma^1,\,\boldsymbol
\gamma^2,\,\boldsymbol\gamma^3$ and through the other four operators 
$\boldsymbol\gamma_0,\,\boldsymbol\gamma_1,\,\boldsymbol\gamma_2,
\,\boldsymbol\gamma_3$:
$$
\hskip -2em
\aligned
\bold H=\frac{i}{24}\sum^3_{p=0}\sum^3_{q=0}\sum^3_{k=0}\sum^3_{m=0}
\omega_{p\kern 0.5pt qkm}\,\boldsymbol\gamma^p\,\boldsymbol\gamma^q
\,\boldsymbol\gamma^k\,\boldsymbol\gamma^m,\\
\bold H=\frac{i}{24}\sum^3_{p=0}\sum^3_{q=0}\sum^3_{k=0}\sum^3_{m=0}
\omega^{p\kern 0.5pt qkm}\,\boldsymbol\gamma_p\,\boldsymbol\gamma_q
\,\boldsymbol\gamma_k\,\boldsymbol\gamma_m.
\endaligned
\mytag{3.7}
$$
Here $24=4\,!$\,. The formulas \mythetag{3.7} are easily proved in special 
frames by means of the formulas \mythetag{2.7}, \mythetag{2.16}, and 
\mythetag{3.6}. In this case they are reduced to 
$$
\xalignat 2
&\hskip -2em
\bold H=i\,\boldsymbol\gamma^0\,\boldsymbol\gamma^1\,\boldsymbol\gamma^2
\,\boldsymbol\gamma^3,
&&\bold H=-i\,\boldsymbol\gamma_0\,\boldsymbol\gamma_1\,\boldsymbol\gamma_2
\,\boldsymbol\gamma_3.
\mytag{3.8}
\endxalignat
$$
By means of the direct calculations we find that
$$
\hskip -2em
i\,\gamma^0\,\gamma^1\,\gamma^2\,\gamma^3=
\Vmatrix 1&0&0&0\\0&1&0&0\\0&0&-1&0\\0&0&0&-1\endVmatrix
=-i\,\gamma_0\,\gamma_1\,\gamma_2\,\gamma_3.
\mytag{3.9}
$$
Comparing \mythetag{3.9} with \mythetag{2.9}, we prove the formulas 
\mythetag{3.8}. Note that \mythetag{3.7} are proper tensorial formulas. 
Having been proved in special frames, they remain valid in an arbitrary 
pair of frames.\par
      Note that in physical literature the operator $\boldsymbol\gamma^5
=-i\,\boldsymbol\gamma^0\,\boldsymbol\gamma^1\,\boldsymbol\gamma^2
\,\boldsymbol\gamma^3$ is introduced (see \S\,22 in \mycite{1}). This 
is another bad tradition since $\boldsymbol\gamma^5$ is not a part of 
the Dirac $\gamma$-field. It is a separate spin-tensorial field 
$\boldsymbol\gamma^5=-\bold H$. I prefer to use the chirality operator 
$\bold H$ instead of the operator $\boldsymbol\gamma^5$.\par
     The chirality operator $\bold H$ anticommutes with the operators 
$\boldsymbol\gamma^0,\,\boldsymbol\gamma^1,\,\boldsymbol\gamma^2,
\,\boldsymbol\gamma^3$ and with the operators $\boldsymbol\gamma_0,
\,\boldsymbol\gamma_1,\,\boldsymbol\gamma_2,\,\boldsymbol\gamma_3$,
i\.\,e\. we have the equalities:
$$
\xalignat 2
&\hskip -2em
\{\bold H,\boldsymbol\gamma^k\}=0,
&&\{\bold H,\boldsymbol\gamma_k\}=0.
\mytag{3.10}
\endxalignat
$$
The equalities \mythetag{3.10} are derived in a special frame by 
means of the formulas \mythetag{2.9}, \mythetag{2.16}, and \mythetag{3.6}.
Then they are extended to arbitrary frame pairs by linearity.\par
     Taking pairs of $\gamma$-operators, we can write the following 
commutation relationships for the operators $\boldsymbol\gamma^0,
\,\boldsymbol\gamma^1,\,\boldsymbol\gamma^2,\,\boldsymbol\gamma^3$ 
and for the operators $\boldsymbol\gamma_0,\,\boldsymbol\gamma_1,
\,\boldsymbol\gamma_2,\,\boldsymbol\gamma_3$:
$$
\xalignat 2
&\hskip -2em
[\kern 0.2pt\bold H,\boldsymbol\gamma^k\,\boldsymbol\gamma^q]=0,
&&[\kern 0.2pt\bold H,\boldsymbol\gamma_k\,\boldsymbol\gamma_q]=0.
\mytag{3.11}
\endxalignat
$$
In the case of three operators we again have the anticommutation
relationships
$$
\xalignat 2
&\hskip -2em
\{\bold H,\boldsymbol\gamma^p\,\boldsymbol\gamma^k
\,\boldsymbol\gamma^q\}=0,
&&\{\bold H,\boldsymbol\gamma_p\,\boldsymbol\gamma_k
\,\boldsymbol\gamma_q\}=0.
\mytag{3.12}
\endxalignat
$$
The formulas \mythetag{3.11} and  \mythetag{3.12} are easily derived 
from the formula \mythetag{3.10}.\par
     Apart from \mythetag{3.11} and  \mythetag{3.12} we need some
additional formulas --- not for commutators and anticommutators, but
for the products of $\gamma$-operators and the chirality operator 
$\bold H$. In the case of two $\gamma$-operators we have
$$
\hskip -2em
\aligned
&\bold H\,\boldsymbol\gamma^p\,\boldsymbol\gamma^q=\bold H
\,g^{p\kern 0.5pt q}-\frac{i}{2}\sum^3_{r=0}\sum^3_{s=0}
\boldsymbol\gamma_r\,\boldsymbol\gamma_s
\ \omega^{rsp\kern 0.5pt q},\\
&\bold H\,\boldsymbol\gamma_p\,\boldsymbol\gamma_q=\bold H
\,g_{p\kern 0.5pt q}-\frac{i}{2}\sum^3_{r=0}\sum^3_{s=0}
\boldsymbol\gamma^r\,\boldsymbol\gamma^s
\ \omega_{rsp\kern 0.5pt q}.
\endaligned
\mytag{3.13}
$$
In the case of three $\gamma$-operators we have a little bit more
complicated formulas:
$$
\hskip -2em
\aligned
&\boldsymbol\gamma^p\,\boldsymbol\gamma^q\,\boldsymbol\gamma^r
=g^{p\kern 0.5pt q}\,\boldsymbol\gamma^r+g^{q\kern 0.2pt r}
\,\boldsymbol\gamma^p-g^{p\kern 0.5pt r}\,\boldsymbol\gamma^q
+i\sum^3_{s=0}\omega^{p\kern 0.5pt qrs}
\,\bold H\,\boldsymbol\gamma_s,\\
&\boldsymbol\gamma_p\,\boldsymbol\gamma_q\,\boldsymbol\gamma_r
=g_{p\kern 0.5pt q}\,\boldsymbol\gamma_r+g_{q\kern 0.2pt r}
\,\boldsymbol\gamma_p-g_{p\kern 0.5pt r}\,\boldsymbol\gamma_q
+i\sum^3_{s=0}\omega_{p\kern 0.5pt qrs}
\,\bold H\,\boldsymbol\gamma^s.
\endaligned
\mytag{3.14}
$$
The formulas \mythetag{3.13} and \mythetag{3.14} are proved by 
choosing some special pair of frames where the operators $\bold H$,
$\boldsymbol\gamma^0,\,\boldsymbol\gamma^1,\,\boldsymbol\gamma^2,
\,\boldsymbol\gamma^3$ and $\boldsymbol\gamma_0,\,\boldsymbol
\gamma_1,\,\boldsymbol\gamma_2,\,\boldsymbol\gamma_3$ are given
by the formulas \mythetag{2.9}, \mythetag{2.16}, \mythetag{3.6},
while the metric tensor $\bold g$ is given by the matrix
\mythetag{2.1}. For the sake of completeness let's write the
following four formulas:
$$
\gather
\hskip -2em
\aligned
&\boldsymbol\gamma^p\,\boldsymbol\gamma^q=\bold 1
\,g^{p\kern 0.5pt q}-\frac{i}{2}\sum^3_{r=0}\sum^3_{s=0}
\bold H\,\boldsymbol\gamma_r\,\boldsymbol\gamma_s
\ \omega^{rsp\kern 0.5pt q},\\
&\boldsymbol\gamma_p\,\boldsymbol\gamma_q=\bold 1
\,g_{p\kern 0.5pt q}-\frac{i}{2}\sum^3_{r=0}\sum^3_{s=0}
\bold H\,\boldsymbol\gamma^r\,\boldsymbol\gamma^s
\ \omega_{rsp\kern 0.5pt q}.
\endaligned
\mytag{3.15}\\
\hskip -2em
\aligned
&\bold H\,\boldsymbol\gamma^p\,\boldsymbol\gamma^q\,\boldsymbol
\gamma^r=g^{p\kern 0.5pt q}\,\bold H\,\boldsymbol\gamma^r
+g^{q\kern 0.2pt r}\,\bold H\,\boldsymbol\gamma^p
-g^{p\kern 0.5pt r}\,\bold H\,\boldsymbol\gamma^q
+i\sum^3_{s=0}\omega^{p\kern 0.5pt qrs}
\,\boldsymbol\gamma_s,\\
&\bold H\,\boldsymbol\gamma_p\,\boldsymbol\gamma_q\,\boldsymbol
\gamma_r=g_{p\kern 0.5pt q}\,\bold H\,\boldsymbol\gamma_r
+g_{q\kern 0.2pt r}\,\bold H\,\boldsymbol\gamma_p
-g_{p\kern 0.5pt r}\,\bold H\,\boldsymbol\gamma_q
+i\sum^3_{s=0}\omega_{p\kern 0.5pt qrs}
\,\boldsymbol\gamma^s.
\endaligned
\mytag{3.16}
\endgather
$$
We derive \mythetag{3.15} and \mythetag{3.16} multiplying both
sides of \mythetag{3.13} and \mythetag{3.14} by $\bold H$ and
taking into account the first identity \mythetag{3.1}.\par
\head
4. Trace formulas.
\endhead
     Note that the traces of all of the $\gamma$-operators 
$\boldsymbol\gamma^0,\,\boldsymbol\gamma^1,\,\boldsymbol\gamma^2,
\,\boldsymbol\gamma^3$ are equal to zero. The same is true for
the $\gamma$-operators with lower spatial index $\boldsymbol
\gamma_0,\,\boldsymbol\gamma_1,\,\boldsymbol\gamma_2,\,\boldsymbol
\gamma_3$ as well as for the products of $\gamma$-operators and
the chirality operator $\bold H$:
$$
\xalignat 2
&\hskip -2em
\tr\boldsymbol\gamma^k=\sum^4_{a=1}\gamma^{ak}_a=0,
&&\tr\boldsymbol\gamma_k=\sum^4_{a=1}\gamma^a_{ak}=0.
\mytag{4.1}\\
\vspace{2ex}
&\hskip -2em
\tr(\bold H\,\boldsymbol\gamma^k)=0,
&&\tr(\bold H\,\boldsymbol\gamma_k)=0.
\mytag{4.2}
\endxalignat
$$
The formulas \mythetag{4.1} and \mythetag{4.2} are proved by direct
calculations with the use of the formulas \mythetag{2.9}, \mythetag{2.16}, 
and \mythetag{3.6}.\par
     Note that $\tr(\bold A\,\bold B)=\tr(\bold B\,\bold A)$. Therefore, for
the traces of the double products of $\gamma$-operators we have the following
formulas:
$$
\xalignat 3
&\tr(\boldsymbol\gamma^p\,\boldsymbol\gamma^q)=4\,g^{p\kern 0.5pt q},
&&\tr(\boldsymbol\gamma_p\,\boldsymbol\gamma_q)=4\,g_{p\kern 0.5pt q},
&&\tr(\boldsymbol\gamma^p\,\boldsymbol\gamma_q)=4\,\delta^p_q.\qquad
\mytag{4.3}
\endxalignat
$$
The formulas \mythetag{4.3} are derived from the anticommutation relationships 
\mythetag{3.3} and \mythetag{3.5}. Applying the formulas \mythetag{4.3} to 
\mythetag{3.13}, we derive
$$
\xalignat 3
&\tr(\bold H\,\boldsymbol\gamma^p\,\boldsymbol\gamma^q)=0,
&&\tr(\bold H\,\boldsymbol\gamma_p\,\boldsymbol\gamma_q)=0,
&&\tr(\bold H\,\boldsymbol\gamma^p\,\boldsymbol\gamma_q)=0.
\qquad
\mytag{4.4}
\endxalignat
$$
In order to calculate the traces of triple products of $\gamma$-operators
we use the formulas \mythetag{3.14}. These formulas immediately yield 
$$
\xalignat 2
&\hskip -2em
\tr(\boldsymbol\gamma^p\,\boldsymbol\gamma^q\,\boldsymbol\gamma^r)=0,
&&\tr(\boldsymbol\gamma_p\,\boldsymbol\gamma_q\,\boldsymbol\gamma_r)=0,
\mytag{4.5}\\
\vspace{1ex}
&\hskip -2em
\tr(\bold H\,\boldsymbol\gamma^p\,\boldsymbol\gamma^q
\,\boldsymbol\gamma^r)=0,
&&\tr(\bold H\,\boldsymbol\gamma_p\,\boldsymbol\gamma_q
\,\boldsymbol\gamma_r)=0.
\mytag{4.6}
\endxalignat
$$
The formulas \mythetag{4.6} are derived from \mythetag{3.16}. In 
deriving both \mythetag{4.5} and \mythetag{4.6} we use the formulas
\mythetag{4.1} and \mythetag{4.2}.\par
     Now let's proceed to the quadruple products of $\gamma$-operators.
For the beginning let's lower the index $r$ in the first formula 
\mythetag{3.14}:
$$
\hskip -2em
\boldsymbol\gamma^p\,\boldsymbol\gamma^q\,\boldsymbol\gamma_r
=g^{p\kern 0.5pt q}\,\boldsymbol\gamma_r+\delta^q_r
\,\boldsymbol\gamma^p-\delta^p_r\,\boldsymbol\gamma^q
+i\sum^3_{n=0}\omega^{p\kern 0.5pt qmn}\,g_{mr}
\,\bold H\,\boldsymbol\gamma_n,
\mytag{4.7}
$$
Then we multiply the equality \mythetag{4.7} on the right by 
$\boldsymbol\gamma_s$:
$$
\boldsymbol\gamma^p\,\boldsymbol\gamma^q\,\boldsymbol\gamma_r
\,\boldsymbol\gamma_s=g^{p\kern 0.5pt q}\,\boldsymbol\gamma_r
\,\boldsymbol\gamma_s+\delta^q_r\,\boldsymbol\gamma^p
\,\boldsymbol\gamma_s-\delta^p_r\,\boldsymbol\gamma^q
\,\boldsymbol\gamma_s+i\sum^3_{n=0}\omega^{p\kern 0.5pt qmn}\,g_{mr}
\,\bold H\,\boldsymbol\gamma_n\,\boldsymbol\gamma_s,
$$
Passing to the traces of both sides of this equality, we take into account
the formulas \mythetag{4.3} and \mythetag{4.4}. As a result we derive 
$$
\hskip -2em
\tr(\boldsymbol\gamma^p\,\boldsymbol\gamma^q\,\boldsymbol\gamma_r
\,\boldsymbol\gamma_s)=4\,g^{p\kern 0.5pt q}\,g_{rs}
+4\,\delta^q_r\,\delta^p_s-4\,\delta^p_r\,\delta^q_s.
\mytag{4.8}
$$
In addition to the formula \mythetag{4.8}, there is a formula for 
$\tr(\bold H\,\boldsymbol\gamma^p\,\boldsymbol\gamma^q\,\boldsymbol\gamma_r
\,\boldsymbol\gamma_s)$. However, in this paper we do not need it.
\head
5. The inverse conversion procedure.
\endhead
     Let's return back to the conversion formula \mythetag{1.1}. Omitting
the spinor indices $a$ and $b$, we can write it as an operator equality:
$$
\hskip -2em
\bold F=u\,\bold 1+v\,\bold H+\sum^3_{k=0}\boldsymbol\gamma^k\,u_k
+\sum^3_{k=0}\bold H\,\boldsymbol\gamma^k\,v_k+\sum^3_{p=0}
\sum^3_{q=0}\boldsymbol\gamma^p\,\boldsymbol\gamma^q\,w_{p\kern 0.4pt q}.
\mytag{5.1}
$$
Applying the formulas \mythetag{4.1}, \mythetag{4.2}, and \mythetag{4.3}
to \mythetag{5.1}, taking into account that 
$$
\xalignat 2
&\hskip -2em
\tr\bold 1=4,
&&\tr\bold H=0,
\mytag{5.2}
\endxalignat 
$$
and remembering the skew-symmetry of $w_{p\kern 0.4pt q}$, we derive
$$
\hskip -2em
u=\frac{1}{4}\tr\bold F.
\mytag{5.3}
$$
The inverse conversion procedure is a series of formulas expressing
$u$, $v$, $u_k$, $v_k$, and $w_{p\kern 0.4pt q}$ through $\bold F$.
The formula \mythetag{5.3} is the first formula in such a series.
Here is the second formula. It expresses $v$ through $\bold F$:
$$
\hskip -2em
v=\frac{1}{4}\tr(\bold H\,\bold F).
\mytag{5.4}
$$
The formula \mythetag{5.4} is derived from \mythetag{5.1} with the use
of the formulas  \mythetag{5.2}, \mythetag{4.1},  \mythetag{4.2},
\mythetag{4.4}, and \mythetag{3.1}.\par
     The third conversion formula should express the components of the
covector field $\bold u$ through $\bold F$. We derive it multiplying 
\mythetag{5.1} on the left by $\gamma_k$:
$$
\hskip -2em
u_k=\frac{1}{4}\tr(\boldsymbol\gamma_k\,\bold F).
\mytag{5.5}
$$
In deriving \mythetag{5.5} we use the formulas \mythetag{5.2},
\mythetag{3.10}, \mythetag{4.1},  \mythetag{4.2}, \mythetag{4.3},
\mythetag{4.4}, \mythetag{4.5}, and \mythetag{3.1}. The fourth 
conversion formula is similar to \mythetag{5.5}:
$$
\hskip -2em
v_k=\frac{1}{4}\tr(\boldsymbol\gamma_k\,\bold H\,\bold F).
\mytag{5.6}
$$
In order to derive the fifth conversion formula we multiply 
the formula \mythetag{5.1} on the left by $\boldsymbol\gamma_q
\,\boldsymbol\gamma_p$. Taking the traces of both sides, then 
we get
$$
\hskip -2em
w_{p\kern 0.4pt q}=\frac{1}{16}\tr(\boldsymbol\gamma_q\,\boldsymbol
\gamma_p\,\bold F)-\frac{1}{16}\tr(\boldsymbol\gamma_p\,\boldsymbol
\gamma_q\,\bold F).
\mytag{5.7}
$$
In deriving \mythetag{5.7} we use the formulas \mythetag{4.6},
\mythetag{4.7} and take into account the skew symmetry of 
$w_{p\kern 0.4pt q}$ with respect to the indices $p$ and $q$.
\par
     The formulas \mythetag{5.3}, \mythetag{5.4}, \mythetag{5.5}, 
\mythetag{5.6}, and \mythetag{5.7} constitute the inverse conversion
procedure. They prove that the mapping \mythetag{5.1}, which produces
a spin-operator $\bold F$ from a collection of purely spatial fields
$u,\,v,\,\bold u,\,\bold v,\bold w$, is bijective.
\head
6. Symmetric and skew-symmetric operators.
\endhead
    Note that the spinor metric $\bold d$ given by the matrices 
\mythetag{2.8} in orthonormal spinor frames defines a skew-symmetric
bilinear form in fibers of the Dirac bundle $DM$:
$$
d(\boldsymbol\psi,\boldsymbol\phi)=\sum^4_{a=1}\sum^4_{b=1}
d_{ab}\,\psi^a\,\phi^b.
\mytag{6.1}
$$ 
\mydefinition{6.1} A spin-operator field $\bold F$ is called a 
\pagebreak {\it symmetric operator\/} if it is symmetric with respect 
to the bilinear form \mythetag{6.1}, i\.\,e\. if $d(\bold F\boldsymbol
\psi,\boldsymbol\phi)=d(\boldsymbol\psi,\bold F\boldsymbol\phi)$ 
for any two spinor fields $\boldsymbol\psi$ and $\boldsymbol\phi$.
\enddefinition
\mydefinition{6.2} A spin-operator field $\bold F$ is called a 
{\it skew-symmetric operator\/} if it is symmetric with respect 
to the bilinear form \mythetag{6.1}, i\.\,e\. if $d(\bold F\boldsymbol
\psi,\boldsymbol\phi)=-d(\boldsymbol\psi,\bold F\boldsymbol\phi)$ 
for any two spinor fields $\boldsymbol\psi$ and $\boldsymbol\phi$.
\enddefinition
\noindent
In the coordinate form the symmetry and skew-symmetry conditions are written 
as
$$
\xalignat 2
&\sum^4_{c=1}F^c_a\,d_{cb}=\sum^4_{c=1}d_{ac}\,F^c_b,
&&\sum^4_{c=1}F^c_a\,d_{cb}=-\sum^4_{c=1}d_{ac}\,F^c_b.
\qquad
\mytag{6.2}
\endxalignat
$$
Using \mythetag{6.2}, we easily prove that the unit operator $\bold 1$ and 
the chirality operator $\bold H$ are symmetric, while the $\gamma$-operators
$\boldsymbol\gamma^k$ and $\boldsymbol\gamma_k$ are skew-symmetric. The
products $\bold H\,\boldsymbol\gamma^k$ and $\bold H\,\boldsymbol\gamma_k$
are symmetric. As for the products $\boldsymbol\gamma^p\,\boldsymbol\gamma^q$
and $\boldsymbol\gamma_p\,\boldsymbol\gamma_q$, they have both symmetric and
skew-symmetric components:
$$
\xalignat 2
&(\boldsymbol\gamma^p\,\boldsymbol\gamma^q)_{\sssize\text{sym}}
=\bold 1\,g^{p\kern 0.4pt q},
&&(\boldsymbol\gamma^p\,\boldsymbol\gamma^q)_{\sssize\text{skew}}
=-\frac{i}{2}\sum^3_{r=0}\sum^3_{s=0}
\bold H\,\boldsymbol\gamma_r\,\boldsymbol\gamma_s
\ \omega^{rsp\kern 0.5pt q},\qquad\\
\vspace{-1.5ex}
\mytag{6.3}\\
\vspace{-1.5ex}
&(\boldsymbol\gamma_p\,\boldsymbol\gamma_q)_{\sssize\text{sym}}
=\bold 1\,g_{p\kern 0.4pt q},
&&(\boldsymbol\gamma_p\,\boldsymbol\gamma_q)_{\sssize\text{skew}}
=-\frac{i}{2}\sum^3_{r=0}\sum^3_{s=0}
\bold H\,\boldsymbol\gamma^r\,\boldsymbol\gamma^s
\ \omega_{rsp\kern 0.5pt q}.\qquad
\endxalignat
$$
The formula \mythetag{6.3} is easily derived from \mythetag{3.15}. Note that
the symmetric parts of the products $\boldsymbol\gamma^p\,\boldsymbol\gamma^q$
and $\boldsymbol\gamma_p\,\boldsymbol\gamma_q$ are symmetric with respect to 
the indices $p$ and $q$, while their skew-symmetric parts are skew-symmetric 
with respect to these indices. This is not a general rule, but a pure 
coincidence in this particular case since the operator symmetrization and the 
operator alternation for spin-operators are not the same as the symmetrization 
and alternation for spatial indices.\par 
\mytheorem{6.1} The spin-operator $\bold F$ given by the formula \mythetag{5.1}
is a symmetric operator if and only if $u_k=0$ and $w_{p\kern 0.4pt q}=0$.
\endproclaim 
\mytheorem{6.2}The spin-operator $\bold F$ given by the formula \mythetag{5.1}
is a skew-symmetric operator if and only if $u=0$, $v=0$, and $v_k=0$.
\endproclaim
\noindent
These two theorems are easily proved on the base of the above results in this
section concerning the operators in the right hand side of the expansion 
\mythetag{5.1}.
\head
7. Hermitian and anti-Hermitian operators.
\endhead
     Let's recall that the bundle of Dirac spinors $DM$ is equipped with the
Dirac form $\bold D$. Its components are given by the matrix \mythetag{2.10}
in self-adjoint spinor frames. Using $\bold D$, we can define a sesquilinear
form in fibers of the bundle $DM$:
$$
\hskip -2em
D(\boldsymbol\psi,\boldsymbol\phi)=\sum^4_{a=1}\sum^4_{\bar b=1}D_{a\bar b}
\,\overline{\psi^{\bar b}}\,\phi^a.
\mytag{7.1}
$$
The sesquilinear form \mythetag{7.1} is not positive. Its signature is 
$(+,+,-,-)$.\par
\mydefinition{7.1} A spin-operator field $\bold F$ is called a 
{\it Hermitian operator\/} if it is Hermitian with respect 
to the sesquilinear form \mythetag{7.1}, i\.\,e\. if $D(\bold F\boldsymbol
\psi,\boldsymbol\phi)=D(\boldsymbol\psi,\bold F\boldsymbol\phi)$ 
for any two spinor fields $\boldsymbol\psi$ and $\boldsymbol\phi$.
\enddefinition
\mydefinition{7.2} A spin-operator field $\bold F$ is called a 
{\it anti-Hermitian operator\/} if it is anti-Hermitian  with respect 
to the sesquilinear form \mythetag{7.1}, i\.\,e\. if for for any two spinor 
fields $\boldsymbol\psi$ and $\boldsymbol\phi$ we have $D(\bold F\boldsymbol
\psi,\boldsymbol\phi)=-D(\boldsymbol\psi,\bold F\boldsymbol\phi)$.
\enddefinition
In the coordinate form the conditions of being Hermitian and anti-Hermitian 
for a spin-operator field $\bold F$ are written as follows:
$$
\xalignat 2
&\sum^4_{\bar c=1}\raise 2.2pt \hbox{$\overline{\lower 2pt
\hbox{$F$}^{\bar c}_{\bar b}}$}\,D_{a\bar c}=
\sum^4_{c=1}D_{c\bar b}\,F^c_a,
&&\sum^4_{\bar c=1}\raise 2.2pt \hbox{$\overline{\lower 2pt
\hbox{$F$}^{\bar c}_{\bar b}}$}\,D_{a\bar c}=
-\sum^4_{c=1}D_{c\bar b}\,F^c_a.\qquad
\mytag{7.2}
\endxalignat
$$
Using \mythetag{7.2}, we easily prove that the unit operator $\bold 1$
is Hermitian, the chirality operator $\bold H$ is anti-Hermitian, and
the $\gamma$-operators $\boldsymbol\gamma^k$ and $\boldsymbol\gamma_k$
are Hermitian operators. The products $\bold H\boldsymbol\gamma^k$ and
$\bold H\boldsymbol\gamma_k$ are also Hermitian operators. The double
products $\boldsymbol\gamma^p\,\boldsymbol\gamma^q$ and $\boldsymbol
\gamma_p\,\boldsymbol\gamma_q$ have both Hermitian and anti-Hermitian
parts:
$$
\xalignat 2
&(\boldsymbol\gamma^p\,\boldsymbol\gamma^q)_{\sssize\text{Herm}}
=\bold 1\,g^{p\kern 0.4pt q},
&&(\boldsymbol\gamma^p\,\boldsymbol\gamma^q)_{\sssize\text{anti}}
=-\frac{i}{2}\sum^3_{r=0}\sum^3_{s=0}
\bold H\,\boldsymbol\gamma_r\,\boldsymbol\gamma_s
\ \omega^{rsp\kern 0.5pt q},\qquad\\
\vspace{-1.5ex}
\mytag{7.3}\\
\vspace{-1.5ex}
&(\boldsymbol\gamma_p\,\boldsymbol\gamma_q)_{\sssize\text{Herm}}
=\bold 1\,g_{p\kern 0.4pt q},
&&(\boldsymbol\gamma_p\,\boldsymbol\gamma_q)_{\sssize\text{anti}}
=-\frac{i}{2}\sum^3_{r=0}\sum^3_{s=0}
\bold H\,\boldsymbol\gamma^r\,\boldsymbol\gamma^s
\ \omega_{rsp\kern 0.5pt q}.\qquad
\endxalignat
$$
Comparing \mythetag{7.3} with \mythetag{6.3}, we see that the subdivision
of $\boldsymbol\gamma^p\,\boldsymbol\gamma^q$ and $\boldsymbol\gamma_p
\,\boldsymbol\gamma_q$ into Hermitian and anti-Hermitian components does
coincide with their subdivision into symmetric and skew-symmetric parts.
This is not a general rule again, but a pure coincidence in our particular
case.
\mytheorem{7.1} The spin-operator $\bold F$ given by the formula \mythetag{5.1}
is a Hermitian operator if and only if we have 
$$
\xalignat 5
&u=\overline{u},
&&v=-\overline{v},
&&u_k=\overline{u_k},
&&v_k=\overline{v_k},
&&w_{p\kern 0.4pt q}=-\overline{w_{p\kern 0.4pt q}}.
\endxalignat
$$
\endproclaim 
\mytheorem{7.2}The spin-operator $\bold F$ given by the formula \mythetag{5.1}
is an anti-Hermi\-tian operator if and only if we have
$$
\xalignat 5
&u=-\overline{u},
&&v=\overline{v},
&&u_k=-\overline{u_k},
&&v_k=-\overline{v_k},
&&w_{p\kern 0.4pt q}=\overline{w_{p\kern 0.4pt q}}.
\endxalignat
$$
\endproclaim
The proof of these two theorems is obvious. Indeed, multiplying a Hermitian 
and an anti-Hermitian operators by a real scalar, we again get a Hermitian 
operator and an anti-Hermitian operator respectively. Multiplying these 
operators by an imaginary scalar, we convert a Hermitian operator into an 
anti-Hermitian operator and vice versa. Therefore, in order to get a Hermitian
operator $\bold F$ the coefficients of Hermitian operators in the expansion
\mythetag{5.1} should be reals, while the coefficients of anti-Hermitian
operators in \mythetag{5.1} should be imaginary numbers. 
\head
8. Commutator equations.
\endhead
    Assume that some spatial frame $\boldsymbol\Upsilon_0,\,\boldsymbol
\Upsilon_1,\,\boldsymbol\Upsilon_2,\,\boldsymbol\Upsilon_3$ is fixed. 
Then the Dirac $\gamma$-field is subdivided into four $\gamma$-operators 
$\boldsymbol\gamma^0,\,\boldsymbol\gamma^1,\,\boldsymbol\gamma^2,
\,\boldsymbol\gamma^3$ or equivalently $\boldsymbol\gamma_0,\,\boldsymbol
\gamma_1,\,\boldsymbol\gamma_2,\,\boldsymbol\gamma_3$. Under this assumption
we consider the following commutator equations:
$$
\hskip -2em
[\kern 0.2pt\bold F,\boldsymbol\gamma_m]=\bold V_{\!m}.
\mytag{8.1}
$$
Here $\bold F$ is an undetermined operator, while $\bold V_{\!m}$ are 
some given operators that constitute a spin-tensorial field of the type 
$(1,1|0,0|0,1)$.\par
      For the beginning we consider the special case where $\bold V_{\!m}=0$.
In this case the operator $\bold F$ in \mythetag{8.1} should commute with
the operators $\boldsymbol\gamma_0,\,\boldsymbol\gamma_1,\,\boldsymbol
\gamma_2,\,\boldsymbol\gamma_3$:
$$
\hskip -2em
[\kern 0.2pt\bold F,\boldsymbol\gamma_m]=0.
\mytag{8.2}
$$
Note that the $\gamma$-matrices \mythetag{3.6} and their products complemented
with the unit matrix span the space of all $4\times 4$ complex matrices. Hence
\mythetag{8.2} means $\bold F$ commutes with all operators acting in fibers 
of the Dirac bundle. Such an operator is scalar, i\.\,e\. it coincides with
the unit operator up to a scalar factor:
$$
\hskip -2em
\bold F=u\,\bold 1.
\mytag{8.3}
$$
The operator \mythetag{8.3} is the general solution of the equation 
\mythetag{8.2}.
\mytheorem{8.1} In the case where the equations \mythetag{8.1} have 
a solution, this solution is unique up to the additive complement 
of the form \mythetag{8.3}.
\endproclaim
      In order to solve the equations \mythetag{8.1} we substitute the
conversion formula \mythetag{5.1} into \mythetag{8.1} for $\bold F$. As 
a result we obtain the equation
$$
\hskip -2em
\gathered
v\,[\kern 0.2pt\bold H,\boldsymbol\gamma_m]+\sum^3_{k=0}u_k
\,[\kern 0.2pt\boldsymbol\gamma^k,\boldsymbol\gamma_m]+\sum^3_{k=0}v_k
\,[\kern 0.2pt\bold H\,\boldsymbol\gamma^k,\boldsymbol\gamma_m]\,+\\
+\sum^3_{p=0}\sum^3_{q=0}w_{p\kern 0.4pt q}\,[\kern 0.2pt\boldsymbol
\gamma^p\,\boldsymbol\gamma^q,\boldsymbol\gamma_m]=\bold V_{\!m}. 
\endgathered
\mytag{8.4}
$$
In order to transform the first term in \mythetag{8.4} we use 
\mythetag{3.10}. This formula yields
$$
\hskip -2em
[\kern 0.2pt\bold H,\boldsymbol\gamma_m]=2\,\bold H\,\boldsymbol\gamma_m
-\{\kern 0.2pt\bold H,\boldsymbol\gamma_m\}=2\,\bold H\,\boldsymbol\gamma_m
=\sum^3_{k=0}2\,g_{mk}\,\bold H\,\boldsymbol\gamma^k.
\mytag{8.5}
$$
For the second term in \mythetag{8.4} we use the skew symmetry of the 
commutator:
$$
\hskip -2em
\gathered
\sum^3_{k=0}u_k\,[\kern 0.2pt\boldsymbol\gamma^k,\boldsymbol\gamma_m]
=\sum^3_{p=0}\sum^3_{q=0}(\boldsymbol\gamma^p\,\boldsymbol\gamma^q
-\boldsymbol\gamma^q\,\boldsymbol\gamma^p)\,g_{qm}=\\
=\sum^3_{p=0}\sum^3_{q=0}\boldsymbol\gamma^p\,\boldsymbol\gamma^q\,
(u_p\,g_{qm}-u_q\,g_{pm}).
\endgathered
\mytag{8.6}
$$
Transforming the third term in \mythetag{8.4}, we use the formulas 
\mythetag{3.3} and \mythetag{3.10}:
$$
\hskip -2em
\gathered
\sum^3_{k=0}v_k\,[\kern 0.2pt\bold H\,\boldsymbol\gamma^k,
\boldsymbol\gamma_m]=\sum^3_{p=0}\sum^3_{q=0}v_p\,(\bold H
\,\boldsymbol\gamma^p\,\boldsymbol\gamma^q-\boldsymbol\gamma^q\,
\bold H\,\boldsymbol\gamma^p)\,g_{qm}=\\
=\sum^3_{p=0}\sum^3_{q=0}v_p\,\bold H\,(\boldsymbol\gamma^p
\,\boldsymbol\gamma^q+\boldsymbol\gamma^q\,\boldsymbol\gamma^p)
\,g_{qm}=2\,v_m\,\bold H.
\endgathered
\mytag{8.7}
$$
And finally, in order to transform the fourth term in \mythetag{8.4}
we apply the formula \mythetag{3.14}. As a result for this term we
derive 
$$
\gathered
\sum^3_{p=0}\sum^3_{q=0}w_{p\kern 0.4pt q}\,[\kern 0.2pt\boldsymbol
\gamma^p\,\boldsymbol\gamma^q,\boldsymbol\gamma_m]=\sum^3_{p=0}
\sum^3_{q=0}\sum^3_{r=0}w_{p\kern 0.4pt q}\,(\boldsymbol\gamma^p
\,\boldsymbol\gamma^q\,\boldsymbol\gamma^r-\boldsymbol\gamma^r
\,\boldsymbol\gamma^p\,\boldsymbol\gamma^q)\,g_{rm}=\\
=\sum^3_{p=0}\sum^3_{q=0}\sum^3_{r=0}2\,w_{p\kern 0.4pt q}
\,(g^{q\kern 0.2pt r}\,\boldsymbol\gamma^p-g^{p\kern 0.5pt r}
\,\boldsymbol\gamma^q)\,g_{rm}=\sum^3_{k=0}4\,w_{km}\,\boldsymbol
\gamma^k.
\endgathered
\mytag{8.8}
$$\par
     Note that the operators $\bold V_{\!m}$ in \mythetag{8.1} can 
also be expressed in the form of \mythetag{5.1}. In order to distinguish
this expression from the original expression \mythetag{5.1} for the
operator $\bold F$ we set the tilde sign  over the coefficients of it:
$$
\bold V_{\!m}=\tilde u_m\,\bold 1+\tilde v_m\,\bold H+\sum^3_{k=0}
\boldsymbol\gamma^k\,\tilde u_{mk}+\sum^3_{k=0}\bold H\,\boldsymbol
\gamma^k\,\tilde v_{mk}+\sum^3_{p=0}\sum^3_{q=0}\boldsymbol\gamma^p
\,\boldsymbol\gamma^q\,\tilde w_{mp\kern 0.4pt q}.\quad
\mytag{8.9}
$$
Using the formulas \mythetag{8.5}, \mythetag{8.6}, \mythetag{8.7}, 
\mythetag{8.8}, \mythetag{8.9}, we prove the following theorem.
\mytheorem{8.2} The commutator equations \mythetag{8.1} are solvable
if and only if the operators $\bold V_{\!m}$ are presented by the 
formula \mythetag{8.9} where $\tilde u_m=0$; \ $\tilde u_{mk}$
is skew symmetric; \ $\tilde v_{mk}=2\,v\,g_{mk}$ for some scalar $v$; 
\ $w_{mp\kern 0.4pt q}=u_p\,g_{qm}-u_q\,g_{pm}$ for the components $u_p$ 
of some covector.
\endproclaim
The theorems~\mythetheorem{8.1} and \mythetheorem{8.2} are helpful in 
studying Lie derivatives for Dirac spinors. However, this is the subject 
for a separate paper.

\enddocument
\end